\documentclass[12pt,a4paper]{article}
\usepackage{amssymb,graphics,amsmath}
\title{}
\date{}

\begin{document}
\begin{center}
\Large{Seven mutually touching inf{\kern0pt}inite cylinders}\\[1cm]
\end{center}
\begin{center}
\normalsize
S\'andor Boz\'oki$^{1}$ \\
\small{
Laboratory on Engineering and Management Intelligence, \\
Research Group of Operations Research and Decision Systems, \\
Institute for Computer Science and Control, \\
Hungarian Academy of Sciences (MTA SZTAKI) \\
Budapest, Hungary \\
\verb|bozoki.sandor@sztaki.mta.hu|  \\
\verb|http://www.sztaki.mta.hu/~bozoki| \\
Research was supported in part by OTKA grant K77420. } \\[5mm]
\end{center}
\begin{center}
\normalsize
Tsung-Lin Lee  \\
\small{
Department of Applied Mathematics, \\
National Sun Yat-sen University \\
Taiwan ROC \\
\verb|leetsung@math.nsysu.edu.tw| \\
\verb|http://www.math.nsysu.edu.tw/~leetsung| \\
Research was supported in part by NSC grant 102-2115-M-110-009. } \\[5mm]
\end{center}
\begin{center}
\normalsize
Lajos R\'onyai \\
\small{
Informatics Laboratory, \\
Institute for Computer Science and Control,  \\
Hungarian Academy of Sciences (MTA SZTAKI); \\
and  \\
Institute of Mathematics,  \\
Budapest University of Technology and Economics,  \\
Budapest, Hungary \\
\verb|ronyai.lajos@sztaki.mta.hu|   \\
\verb|http://www.sztaki.mta.hu/~ronyai|   \\
Research was supported in part by OTKA grants K77476 and NK105645. }\\[5mm]
\end{center}
\begin{center}
\normalsize{28 May 2014}
\end{center}
\footnotetext[1]{corresponding author}
\newpage
\begin{abstract}
We solve a problem of Littlewood: there exist seven
inf{\kern0pt}inite circular cylinders of unit radius which
mutually touch each other. In fact, we exhibit two such sets of
cylinders. Our approach is algebraic and uses symbolic and
numerical computational techniques. We consider a system of
polynomial equations describing the position of the axes of the
cylinders in the 3 dimensional space. To have the same number of
equations (namely 20) as the number of variables, the angle of the
f{\kern0pt}irst two cylinders is f{\kern0pt}ixed to 90 degrees,
and a small family of direction vectors is left out of
consideration. Homotopy continuation method has been applied to
solve the system. The number of paths is about 121 billion, it is
hopeless to follow them all. However, after checking 80 million
paths, two solutions are found. Their validity, i.e., the
existence of exact real solutions close to the approximate
solutions at hand, was verif{\kern0pt}ied with  the
alphaCertif{\kern0pt}ied method as well as by the interval
Krawczyk method.
\end{abstract}
homotopy method,
certif{\kern0pt}ied solutions,
alpha theory,
interval methods 
\section{Littlewood's problem on seven touching in-f{\kern0pt}inite cylinders}

John Edensor Littlewood (\cite{Littlewood1968}, Problem 7 on p.~20)
proposed that
\begin{center}
\parbox{10cm}{
``\emph{Is it possible in 3-space for seven inf{\kern0pt}inite circular
cylinders of unit radius each to touch all the others?
Seven is the number suggested by constants.}''
}
\end{center}
Two cylinders touch each other if their intersection is
either a point or a line.
Ogilvy's book \cite{Ogilvy1962} also includes Littlewood's problem. \\

Finite versions of the problem are discussed as puzzles by
Gardner and they are well known as
6 touching cigarettes \cite[Figure 54 on page 115]{Gardner1988}
and 7 touching cigarettes \cite[Figure 55 on page 115]{Gardner1988}.
The latter works for a ratio  of
length/diameter greater than $7\sqrt{3}/2$. However, as it is noted by Bezdek
\cite{Bezdek2005}, it is still open whether it is possible to
f{\kern0pt}ind 8 or more touching f{\kern0pt}inite identical
cylinders. An arrangement of 5
touching coins (with a \emph{small} ratio of length/diameter) is
also known \cite[Figure 49 on page 110]{Gardner1988} and this fact
suggests that intermediate ratios of length/diameter
could also be  analyzed.  \\

Bezdek \cite{Bezdek2005} showed that 24 is an upper bound
for the number of mutually touching congruent inf{\kern0pt}inite cylinders.
Ambrus and Bezdek \cite{AmbrusBezdek2008} investigated the
proposal of Kuperberg from the early 1990's that contained
8 congruent inf{\kern0pt}inite cylinders.
It is shown that they do not mutually touch each other,
see \cite[Theorem 1 and Figure 1 on page 1804]{AmbrusBezdek2008}
for details.
Brass, Moser and Pach discuss an arrangement of 6 mutually
touching inf{\kern0pt}inite cylinders \cite[page 98]{BrassMoserPach2005}.
In the paper this lower bound is improved to 7. \\

Hereafter, it is assumed that cylinders are inf{\kern0pt}inite and
congruent, their radius is set to 1.
Two cylinders of unit radius touch each other
if and only if the distance of their axes is 2.
Let $C_i$ and $\ell_i$ denote the $i$-th cylinder
and its axis, respectively. In the paper, $i=1,2,\ldots,7.$
The case of parallel cylinders (lines) is excluded
from our analysis. It is left to the reader to show
that if two cylinders are parallel, then the maximum number
of mutually touching cylinders is four. \\

We intend to apply the well-known formula for the distance of two lines
in $\mathbb R^3$. Let
\[
\ell_i(s) = \mathbf{P}_i + s \, {\mathbf{w}_i}
\]
be a parametric representation of line $\ell_i$ for $i=1,\ldots ,7$.
Here $\mathbf{P}_i \in \mathbb{R}^3$ is a point of $\ell_i$,
$\mathbf{w}_i \in \mathbb{R}^3$ is a direction vector
and $s$ is a real parameter.
If lines $\ell_i$ and $\ell_j$ are skew, then their distance
can be obtained as
\begin{equation} \label{eq:LineLineDistance1}
d(\ell_i,\ell_j) = \frac{|(\overrightarrow{\mathbf{P}_i\mathbf{P}_j})\cdot
({\mathbf{w}_i}\times {\mathbf{w}_j})|}{
||{\mathbf{w}_i}\times {\mathbf{w}_j}||},
\end{equation}
where $\cdot$ denotes dot product,
$\times$ denotes cross product and
$||\,\,||$ denotes the Euclidean norm \cite{GreenspanBenney1997,Weisstein-MathWorld}.
Since the cylinders have unit radius, $d(\ell_i,\ell_j) = 2$
for all $i,j = 1,2,\ldots, 7, \, i \neq j$,
we can write equations (\ref{eq:LineLineDistance1}) as
\begin{equation} \label{eq:LineLineDistance2}
 {|(\overrightarrow{\mathbf{P}_i\mathbf{P}_j})\cdot({\mathbf{w}_i}\times {\mathbf{w}_j})|}^2
 - 4 {||{\mathbf{w}_i}\times {\mathbf{w}_j}||^2} = 0.
\end{equation}
In this form we avoid taking square roots. Let us introduce coordinates:
\[
\mathbf{P}_i = (x_i, y_i, z_i),   \qquad
\mathbf{w}_i = (t_i, u_i, v_i).
\]
Then we have
\begin{eqnarray}
\overrightarrow{\mathbf{P}_i\mathbf{P}_j} &=& ( x_j - x_i , y_j - y_i , z_j - z_i), \label{eq:PiPjxixjyiyjzizj} \\
{\mathbf{w}_i}\times {\mathbf{w}_j} &=& ( u_i v_j-v_i u_j , v_i t_j - t_i v_j , t_i u_j - u_i t_j).
                                                                                    \label{eq:wixwjuiujtitjvivj}
\end{eqnarray}
Now we substitute
(\ref{eq:PiPjxixjyiyjzizj})-(\ref{eq:wixwjuiujtitjvivj}) into (\ref{eq:LineLineDistance2}),
and by using the well-known determinantal form of the triple product,
we obtain the equation
\hspace*{-20mm}
\begin{eqnarray}
\det
\left[
\begin{array}{ccc}
x_j - x_i   &   y_j - y_i   &   z_j-z_i     \\
   t_i      &      u_i      &     v_i       \\
   t_j      &      u_j      &     v_j
\end{array}
\right]^2  
-4
\left(
(u_i v_j - v_i u_j)^2+ \right.   \nonumber \\
+ \left. (v_i t_j - t_i v_j)^2+
(t_i u_j - u_i t_j)^2
\right) = 0.  \label{eq:li-lj-xiyizi-tiuivi}
\end{eqnarray}
This is a polynomial equation of degree 6 in 12 variables.
The polynomial on the left is a linear combination of 84 monomials. \\

We call a line horizontal if it is parallel to the
plane $z=0.$
Any arrangement of seven lines can be translated and rotated
to a position in which one of the lines $(\ell_1)$
is horizontal, with direction vector
$\mathbf{w}_1 = (1,0,0)$, and it goes through the point $\mathbf{P}_1(0,0,-1)$.
It can also be assumed that the touching point of
cylinders $C_1$ and $C_2$ is $(0,0,0)$, that is,
$\ell_2$ goes through the point $\mathbf{P}_2(0,0,1)$.
The direction of $(\ell_2)$ is the only degree
of freedom when the f{\kern0pt}irst two lines are
considered. We shall assume, and this is explained later,
that $(\ell_2)$ will be chosen to be orthogonal to
the f{\kern0pt}irst line. We have so far
\begin{eqnarray}
      x_1 = 0, \, y_1 = 0, \, z_1 = -1,      \qquad  t_1 = 1, \, u_1 = 0, \, v_1 = 0; \label{eq:specificationline1}  \\
      x_2 = 0, \, y_2 = 0, \, z_2 = \quad 1, \qquad  t_2 = 0, \, u_2 = 1, \, v_2 = 0. \label{eq:specificationline2}
\end{eqnarray}
We can make some further simplif{\kern0pt}ications.
We may assume without loss of generality that
 $\ell_i$ $(i = 3,\ldots,7)$ is
not horizontal (otherwise it would be parallel to $\ell_1$ or $\ell_2$),
consequently, it goes through the plane $z=k$ for any $k \in \mathbb{R}$.
Let us choose $k=0$ and set
\begin{eqnarray}
z_i = 0 \quad \textrm{for } i = 3,\ldots,7.     \label{eq:z34567}
\end{eqnarray}
Finally, the normalization
of the direction vector of line $\ell_i$ is chosen to be
$t_i + u_i + v_i = 1$ for $i = 3,\ldots,7.$
This is equivalent to
\begin{equation}
v_i = 1-t_i- u_i, \qquad i = 3,\ldots,7.    \label{eq:v34567}
\end{equation}
The normalization is restrictive, it may rule out some valid
solutions, as it excludes
all nonzero direction vectors fulf{\kern0pt}illing $t_i + u_i + v_i = 0.$
However, our aim is to f{\kern0pt}ind \emph{one} solution
rather than an analysis of all solutions.
At this point we leave it
open whether the excluded direction vectors may produce a valid solution. \\

The distance of $\ell_1$ and $\ell_2$ is guaranteed to be 2 by
(\ref{eq:specificationline1})-(\ref{eq:specificationline2}).
Substitute (\ref{eq:specificationline1}), (\ref{eq:z34567})
and (\ref{eq:v34567}) into (\ref{eq:li-lj-xiyizi-tiuivi})
to set the distance of $\ell_1$ and $\ell_j$ $(3 \leq j \leq 7)$:
\begin{eqnarray}
y_j^2 t_j^2
+2 y_j^2 t_j u_j
-2 y_j^2 t_j
+y_j^2 u_j^2
-2 y_j^2 u_j
+y_j^2
+2 y_j t_j u_j
+2 y_j u_j^2
-2 y_j u_j  \qquad   \nonumber \\
-4 t_j^2
-8 t_j u_j
+8 t_j
-7 u_j^2
+8 u_j
-4
= 0,  \qquad j = 3,\ldots,7. \qquad \qquad \qquad \quad \label{eq:line1-linej}
\end{eqnarray}
Substitute (\ref{eq:specificationline2}), (\ref{eq:z34567})
and (\ref{eq:v34567}) into (\ref{eq:li-lj-xiyizi-tiuivi})
to set the distance of $\ell_2$ and $\ell_j$ $(3 \leq j \leq 7)$:
\begin{eqnarray}
x_j^2 t_j^2
+2 x_j^2 t_j u_j
-2 x_j^2 t_j
+x_j^2 u_j^2
-2 x_j^2 u_j
+x_j^2
-2 x_j t_j u_j
-2 x_j t_j^2
+2 x_j t_j \qquad   \nonumber \\
-4 u_j^2
-8 t_j u_j
+8 t_j
-7 t_j^2
+8 u_j
-4
= 0,  \qquad j = 3,\ldots,7. \qquad \qquad \qquad \quad \label{eq:line2-linej}
\end{eqnarray}

Finally, substitute (\ref{eq:z34567}) and (\ref{eq:v34567}) into (\ref{eq:li-lj-xiyizi-tiuivi})
to set the distance of $\ell_i$ and $\ell_j$ $(3 \leq i < j \leq 7)$: \\
$-4x_i y_i t_i u_i t_j u_j$
$+4 x_i x_j t_i u_i t_j u_j$
$+4 x_i y_j t_i u_i t_j u_j$
$+4 y_i x_j t_i u_i t_j u_j$
$+4 y_i y_j t_i u_i t_j u_j$   \\
$-4 x_j y_j t_i u_i t_j u_j$
$-2 x_i^2 t_i u_i t_j u_j$
$-2 y_i^2 t_i u_i t_j u_j$
$-2 x_j^2 t_i u_i t_j u_j$
$-2 y_j^2 t_i u_i t_j u_j$   \\    
$-4 x_i x_j t_i u_i u_j$
$+4 x_i x_j t_i u_j^2$
$+4 x_i x_j u_i^2 t_j$
$-4 x_i x_j u_i t_j u_j$
$+4 y_i y_j t_i^2 u_j$  \\   
$-4 y_i y_j t_i u_i t_j$
$-4 y_i y_j t_i t_j u_j$
$+4 y_i y_j u_i t_j^2$
$+4 x_i x_j u_i u_j$
$+4 y_i y_j t_i t_j$   \\  
$+x_i^2 t_i^2 u_j^2$
$+x_i^2 u_i^2 t_j^2$
$+y_i^2 t_i^2 u_j^2$
$+y_i^2 u_i^2 t_j^2$
$+x_j^2 t_i^2 u_j^2$
$+x_j^2 u_i^2 t_j^2$
$+y_j^2 t_i^2 u_j^2$
$+y_j^2 u_i^2 t_j^2$  \\    
$+2 x_i y_i t_i^2 u_j^2$
$+2 x_i y_i u_i^2 t_j^2$
$-2 x_i x_j t_i^2 u_j^2$
$-2 x_i x_j u_i^2 t_j^2$
$-2 x_i y_j t_i^2 u_j^2$
$-2 x_i y_j u_i^2 t_j^2$  \\  
$-2 y_i x_j t_i^2 u_j^2$
$-2 y_i x_j u_i^2 t_j^2$
$-2 y_i y_j t_i^2 u_j^2$
$-2 y_i y_j u_i^2 t_j^2$
$+2 x_j y_j t_i^2 u_j^2$
$+2 x_j y_j u_i^2 t_j^2$   \\  
$-2 x_i y_i t_i^2 u_j$
$-2 x_i y_i t_i u_j^2$
$+2 x_i y_j t_i^2 u_j$
$+2 x_i y_j t_i u_j^2$
$+2 x_i y_j u_i^2 t_j$
$+2 x_i y_j u_i t_j^2$    \\  
$-2 x_i y_i u_i^2 t_j$
$-2 x_i y_i u_i t_j^2$
$+2 y_i x_j t_i^2 u_j$
$+2 y_i x_j t_i u_j^2$
$+2 y_i x_j u_i^2 t_j$
$+2 y_i x_j u_i t_j^2$   \\  
$-2 x_j y_j t_i^2 u_j$
$-2 x_j y_j t_i u_j^2$
$-2 x_j y_j u_i^2 t_j$
$-2 x_j y_j u_i t_j^2$
$-2 x_i^2 t_i u_j^2$
$-2 x_i^2 u_i^2 t_j$
$-2 y_i^2 t_i^2 u_j$         \\  
$-2 y_i^2 u_i t_j^2$
$-2 x_j^2 t_i u_j^2$
$-2 x_j^2 u_i^2 t_j$
$-2 y_j^2 t_i^2 u_j$
$-2 y_j^2 u_i t_j^2$
$+2 x_i^2 t_i u_i u_j$
$+2 x_i^2 u_i t_j u_j$      \\  
$+2 y_i^2 t_i u_i t_j$
$+2 y_i^2 t_i t_j u_j$
$+2 x_j^2 t_i u_i u_j$
$+2 x_j^2 u_i t_j u_j$
$+2 y_j^2 t_i u_i t_j$
$+2 y_j^2 t_i t_j u_j$          \\  
$+2 x_i y_i t_i u_i t_j$
$+2 x_i y_i t_i u_i u_j$
$+2 x_i y_i t_i t_j u_j$
$+2 x_i y_i u_i t_j u_j$
$-2 x_i y_j t_i u_i t_j$      \\  
$-2 x_i y_j t_i u_i u_j$
$-2 x_i y_j t_i t_j u_j$
$-2 x_i y_j u_i t_j u_j$
$-2 y_i x_j t_i u_i t_j$
$-2 y_i x_j t_i u_i u_j$       \\  
$-2 y_i x_j t_i t_j u_j$
$-2 y_i x_j u_i t_j u_j$
$+2 x_j y_j t_i u_i t_j$
$+2 x_j y_j t_i u_i u_j$
$+2 x_j y_j t_i t_j u_j$         \\  
$+2 x_j y_j u_i t_j u_j$
$-2 x_i^2 u_i u_j$
$-2 y_i^2 t_i t_j$
$-2 x_j^2 u_i u_j$
$-2 y_j^2 t_i t_j$
$-2 x_i y_i t_i u_i$
$+2 x_i y_i t_i u_j$           \\ 
$+2 x_i y_i u_i t_j$
$-2 x_i y_i t_j u_j$
$+2 x_i y_j t_i u_i$
$-2 x_i y_j t_i u_j$
$-2 x_i y_j u_i t_j$
$+2 x_i y_j t_j u_j$        \\ 
$+2 y_i x_j t_i u_i$
$-2 y_i x_j t_i u_j$
$-2 y_i x_j u_i t_j$
$+2 y_i x_j t_j u_j$
$-2 x_j y_j t_i u_i$
$+2 x_j y_j t_i u_j$          \\ 
$+2 x_j y_j u_i t_j$
$-2 x_j y_j t_j u_j$
$-2 x_i x_j u_i^2$
$-2 x_i x_j u_j^2$          
$-2 y_i y_j t_j^2$
$-2 y_i y_j t_i^2$
$+24 t_i u_i t_j u_j$   \\ 
$+x_i^2 u_i^2$
$+x_i^2 u_j^2$
$+y_i^2 t_i^2$
$+y_i^2 t_j^2$
$+x_j^2 u_i^2$
$+x_j^2 u_j^2$
$+y_j^2 t_i^2$
$+y_j^2 t_j^2$
$-12 t_i^2 u_j^2$
$-12 u_i^2 t_j^2$         \\ 
$-4 t_i^2$
$-4 u_i^2$
$-4 t_j^2$
$-4 u_j^2 $
$-8 t_i u_i t_j$
$-8 t_i u_i u_j$
$-8 t_i t_j u_j$
$+8 t_i u_j^2$
$+8 t_i^2 u_j$
$+8 u_i^2 t_j$       \\
$+8 u_i t_j^2$            
$-8 u_i t_j u_j$
$+8 t_i t_j$
$+8 u_i u_j = 0,  \qquad i = 3,\ldots,6, \, \,  j = i+1,\ldots,7.$                 
\begin{equation} \label{eq:linei-linej}
\end{equation}

As the f{\kern0pt}irst two lines $\ell_1, \ell_2$ are f{\kern0pt}ixed, the aim is to f{\kern0pt}ind
f{\kern0pt}ive lines $\ell_3, \ldots, \ell_7$ such that the distance of each pair of lines is 2.
System (\ref{eq:line1-linej})-(\ref{eq:linei-linej})
has 20 equations and 20 variables ($x_i, y_i, t_i, u_i, $ $i=3,\ldots,7$).
Each equation is a multivariate polynomial equation.
Note that without f{\kern0pt}ixing the angle of lines
$\ell_1$ and $\ell_2$ at a given value, 90 degrees by our choice,
we would have a system of 20 equations and 21 variables and
we would lose the chance of f{\kern0pt}inding isolated roots. \\

The numerical solution is presented in the next section.
We emphasize here that, via methods such as alphaCertif{\kern0pt}ied discussed in
subsection 3.1, the numbers themselves given in Table 1 prove the
existence of real solutions of system
(\ref{eq:line1-linej})-(\ref{eq:linei-linej})
and hence solve Littlewood's problem.
Nevertheless, we also outline the method of computing
the approximate solutions in Table 1.

\section{Solving the polynomial system by the polyhedral homotopy
continuation method}

The polyhedral homotopy continuation method is developed in \cite%
{HuberSturmfels1995} to approximate all isolated zeros of a polynomial
system and is well implemented in software HOM4PS-2.0 \cite{LeeLiTsai2008}.
The numerical experiments show that the method is ef{\kern0pt}f{\kern0pt}icient and reliable.
More importantly, it can handle the large scale polynomial systems such as
the system (\ref{eq:line1-linej})-(\ref{eq:linei-linej}).

For a system of polynomials $P(\mathbf{x})=\left( p_{1}(
\mathbf{x}),\dots ,p_{n}(\mathbf{x})\right) \,$
with $\,\mathbf{x}=\left( x_{1},\dots ,x_{n}\right) ,\,$ write
\begin{equation*}
p_{j}(\mathbf{x})=\underset{\mathbf{a}\in \mathcal{S}_{j}}{\sum }c_{j,
\mathbf{a}}\mathbf{x}^{\mathbf{a}}\text{, \ \ \ \ \
\ \ \ }j=1,\dots ,n,
\end{equation*}%
where $\,\mathbf{a}
=\left( a_{1},\dots ,a_{n}\right) \in (\mathbb{N}\cup \{0\})^{n}$,
$c_{j,\mathbf{a}}\in \mathbb{C}^{\ast }=\mathbb{C}\backslash \{0\}$, $%
\mathbf{x}^{\mathbf{a}}=x_{1}^{a_{1}}\cdots x_{n}^{a_{n}}$, and $
\mathcal{S}_{j}\subset (\mathbb{N}\cup
\{0\})^{n}$
is f{\kern0pt}inite.

Let \thinspace $\omega _{j}:\mathcal{S}_{j}\rightarrow \mathbb{R}$ be a
random lifting function on $\mathcal{S}_{j}$\ which lifts $\mathcal{S}_{j}$\
to its graph $\hat{\mathcal{S}}_{j}=\{\hat{
\mathbf{a}}=(\mathbf{a},\omega _{j}(\mathbf{a}))|
\mathbf{a}\mathbf{%
\in }\mathcal{S}_{j}\}\subset \mathbb{R}^{n+1}$.\thinspace\ A collection of
pairs $(\{
\mathbf{a}_{1},\mathbf{a}^{\prime}_{1}\},\dots ,
\{\mathbf{a}_{n},\mathbf{a}^{\prime}_{n}\}), \text{%
where}\,
\, \mathbf{a}_{j},\mathbf{a}_{j}^{\prime} \in \mathcal{S}_{j},$
is called a
\textit{mixed cell} if there exists $\hat{\boldsymbol\alpha}=(\boldsymbol\alpha ,1)\in \mathbb{R}%
^{n+1}$\ such that
\begin{equation*}
\langle \hat{\mathbf{a}}_{j},\hat{\boldsymbol\alpha }\rangle =\langle \hat{
\mathbf{a}}_{j}^{\prime },\hat{%
\boldsymbol\alpha }\rangle <\langle \hat{
\mathbf{a}},\hat{\boldsymbol\alpha }\rangle \text{\ \ \ for all}%
\;\;
\mathbf{a}\in \mathcal{S}_{j}\backslash \{
\mathbf{a}_{j},
\mathbf{a}_{j}^{\prime }\}\text{, \ }%
j=1,\dots ,n.
\end{equation*}%
Here, $\langle \;,\;\rangle $\thinspace stands for the usual inner product
in the Euclidean space $\mathbb{R}^{n+1}.$ It is well-known that the number of mixed cells of a
polynomial system is f{\kern0pt}inite \cite{HuberSturmfels1995}. Those mixed cells
play an important role in constructing the polyhedral homotopy.

Consider a given mixed cell $\,C=(\{
\mathbf{a}_{11},\mathbf{a}_{12}\},\dots ,\{
\mathbf{a}_{n1},\mathbf{a}_{n2}\})\,$%
with inner normal $\,\boldsymbol\alpha \in \mathbb{R}^{n}$,
where
$\mathbf{a}_{j1},\mathbf{a}_{j2} \in \mathcal{S}_{j}$
for each $\;j=1,\dots ,n$.
Let
$\tilde{c}_{j,\mathbf{a}}$
be a randomly chosen number in $\mathbb{C}$, and denote%
\begin{equation*}
\beta _{j}=\min_{\mathbf{a}\in \mathcal{S}_{j}}\langle\hat{
\mathbf{a}},\hat{\boldsymbol\alpha}\rangle=
\langle\hat{\mathbf{a}}_{j1},%
\hat{\boldsymbol\alpha}\rangle=
\langle\hat{
\mathbf{a}}_{j2},\hat{\boldsymbol\alpha}\rangle.
\end{equation*}%
HOM4PS-2.0 constructs the homotopy to be $H(
\mathbf{x},t)=(h_{1}(\mathbf{x},t),\dots
,h_{n}(\mathbf{x},t))$, $t\in (-\infty ,0]$, where%
\begin{equation*}
h_{j}(\mathbf{x},t)=\sum_{\mathbf{a}
\in \mathcal{S}_{j}}[(1-e^{t})\tilde{c}_{j,\mathbf{a}}+e^{t}c_{j,%
\,
\mathbf{a}}]
\mathbf{x}^{\mathbf{a}}e^{t\cdot (\langle\hat{
\mathbf{a}},\hat{\boldsymbol\alpha}\rangle-\beta _{j})}\quad \mathrm{for}%
~\;j=1,\dots ,n.
\end{equation*}%
Note that $H(\mathbf{x},0)=P(\mathbf{x})$.
When $t$ goes to $-\infty $, $H(
\mathbf{x},t)\;$becomes a
binomial system
\begin{equation*}
\left\{
\begin{array}{c}
\tilde{c}_{11}
\mathbf{x}^{\mathbf{a}_{11}}+\tilde{c}_{12}
\mathbf{x}^{\mathbf{a}_{12}}=0 \\
\vdots  \\
\tilde{c}_{n1}
\mathbf{x}^{\mathbf{a}_{n1}}+\tilde{c}_{n2}
\mathbf{x}^{\mathbf{a}_{n2}}=0%
\end{array}%
\right.
\end{equation*}%
having $\;\left\vert \det \left(
\mathbf{a}_{11}-\mathbf{a}_{12},\dots ,
\mathbf{a}_{n1}-\mathbf{a}_{n2}\right)
\right\vert \,$
nonsingular
 isolated solutions which provide the
starting points for tracking the solution paths of
$\,H(\mathbf{x},t)=0$ from $%
t=-\infty $ to $0$. For the details of the
algorithm for tracking the solution paths, see \cite{LeeLiTsai2008}.

\bigskip

The polynomial system (\ref{eq:line1-linej})-(\ref{eq:linei-linej})
consists of $20$ equations in $20$ variables. We
obtain $180,734$ mixed cells of the system by software MixedVol-2.0 \cite{ChenLeeLi2014}, which
provide $121,098,993,664$ homotopy curves to be tracked. In order to track
so many curves ef{\kern0pt}f{\kern0pt}iciently, we use the subroutines in the TBB library
(Thread Building Blocks) to distribute data over multiple cores for parallel
computation. Employing total $12$ cores in 2 Intel Xeon X5650 2.66 GHz CPUs,
20 million curves are completed in a month.
The f{\kern0pt}irst real solution was
found after tracking 80 million paths,
and the second one was found after tracking
another 25 million paths.

\begin{center}
\begin{tabular}{|l|r@{.}l|r@{.}l|}
\hline
& \multicolumn{2}{c|}{f{\kern0pt}irst solution} & \multicolumn{2}{c|}{second solution}
\\ \hline
$x_3 $ & 11 & 675771704477 & 2 & 075088491891 \\
$y_3 $ & $-4$ & 124414157636 & $-2$ & 036516392124 \\
$t_3 $ & 0 & 704116159640 & $-0$ & 030209763440 \\
$u_3 $ & 0 & 235129952793 & 0 & 599691085438 \\ \hline
$x_4 $ & 3 & 802878122730 & $-2$ & 688893665930 \\
$y_4 $ & $-2$ & 910611127075 & 4 & 070505903499 \\
$t_4 $ & 0 & 895623427074 & 0 & 184499043058 \\
$u_4 $ & $-0$ & 149726023342 & 0 & 426965115851 \\ \hline
$x_5 $ & 8 & 311818491659 & $-4$ & 033142850644 \\
$y_5 $ & $-1$ & 732276613733 & $-2$ & 655943449984 \\
$t_5 $ & 2 & 515897624878 & 0 & 251380280590 \\
$u_5 $ & $-0$ & 566129665502 & 0 & 516678258430 \\ \hline
$x_6 $ & $-6$ & 487945444917 & 6 & 311134419772 \\
$y_6 $ & $-8$ & 537495065091 & $-5$ & 229892181735 \\
$t_6 $ & 0 & 785632006191 & $-0$ & 474742889365 \\
$u_6 $ & 0 & 338461562103 & 1 & 230302197822 \\ \hline
$x_7 $ & $-3$ & 168475045360 & 3 & 914613907006 \\
$y_7 $ & $-2$ & 459640638529 & $-7$ & 881492743224 \\
$t_7 $ & 0 & 192767499267 & 1 & 698198197367 \\
$u_7 $ & 0 & 536724141124 & $-1$ & 164062857743 \\ \hline
\end{tabular}

\textbf{Table 1.} Two solutions of system (\ref{eq:line1-linej})-(\ref%
{eq:linei-linej}) by HOM4PS-2.0
\end{center}

Since system (\ref{eq:line1-linej})-(\ref{eq:linei-linej}) is
symmetric in the f{\kern0pt}ive 4-tuples $(x_j, y_j, t_j, u_j), \, j=3,\ldots,7,$ each
solution represents a family of $5!=120$ solutions, all of them resulting in
the same arrangements of the cylinders. The two solutions in Table 1 are
obviously not permutations of each other. However, due to the rotational and
ref{\kern0pt}lectional symmetries of the orthogonally f{\kern0pt}ixed pair of cylinders $C_1,
C_2 $, any arrangement represents a family of 8 congruent arrangements. In
order to show that the two solutions in Table 1 are non-congruent
arrangements, we have computed the angles between the pairs of cylinders.
The two sets of pairwise angles are disjoint except for the right angle of $%
C_1, C_2$. Consequently, the two arrangements in Figure 1 and 2 are not
congruent.

\newpage
\unitlength 1mm
\begin{center}
\begin{picture}(120,150)
\put(-13,80){\resizebox{140mm}{!}{\rotatebox{0}{
\includegraphics{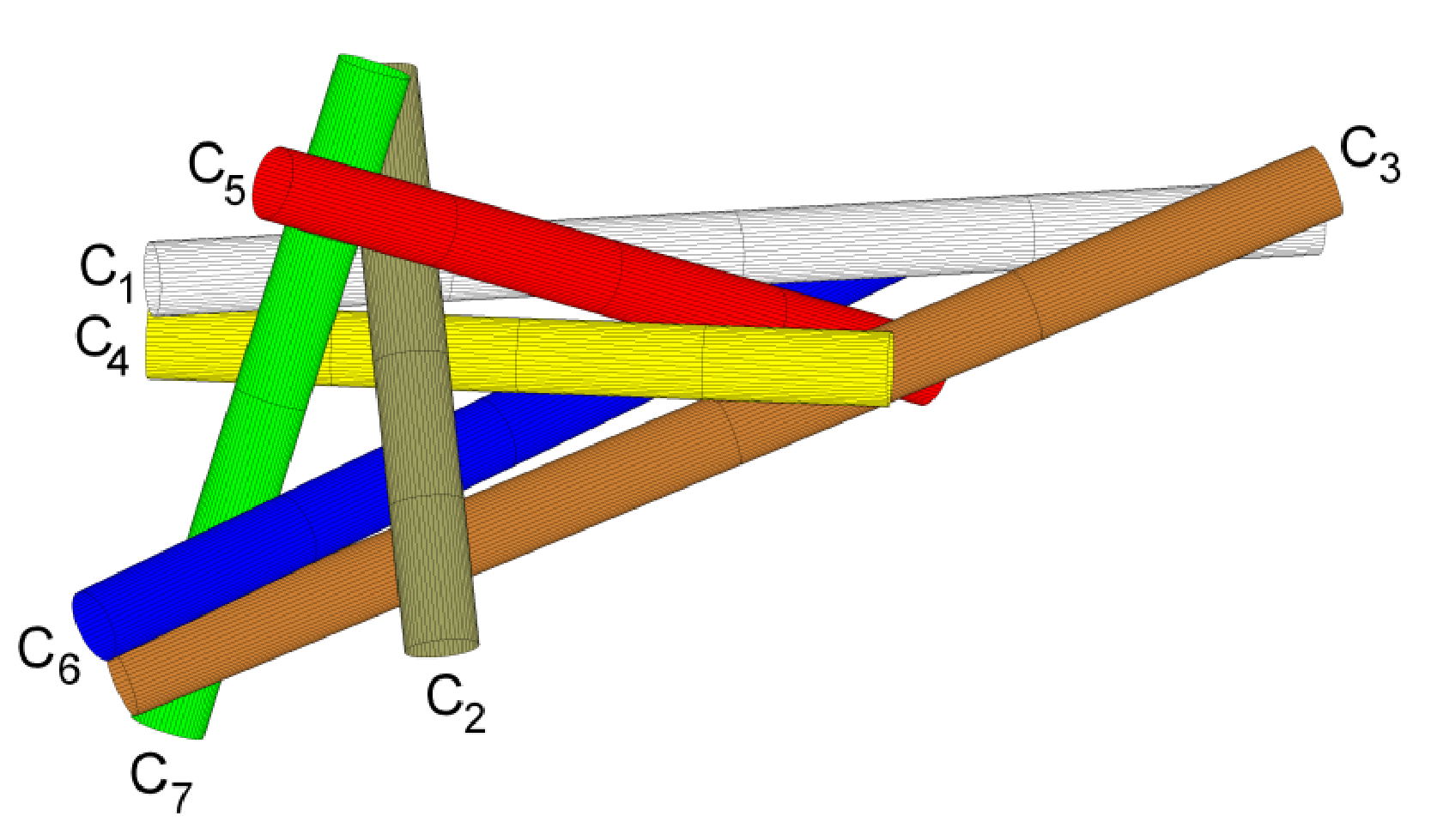}}}}
\put(-10,70){\makebox{\textbf{Figure 1.} The f{\kern0pt}irst set of seven mutually touching inf{\kern0pt}inite cylinders }}
\put(5,-30){\resizebox{100mm}{!}{\rotatebox{0}{
\includegraphics{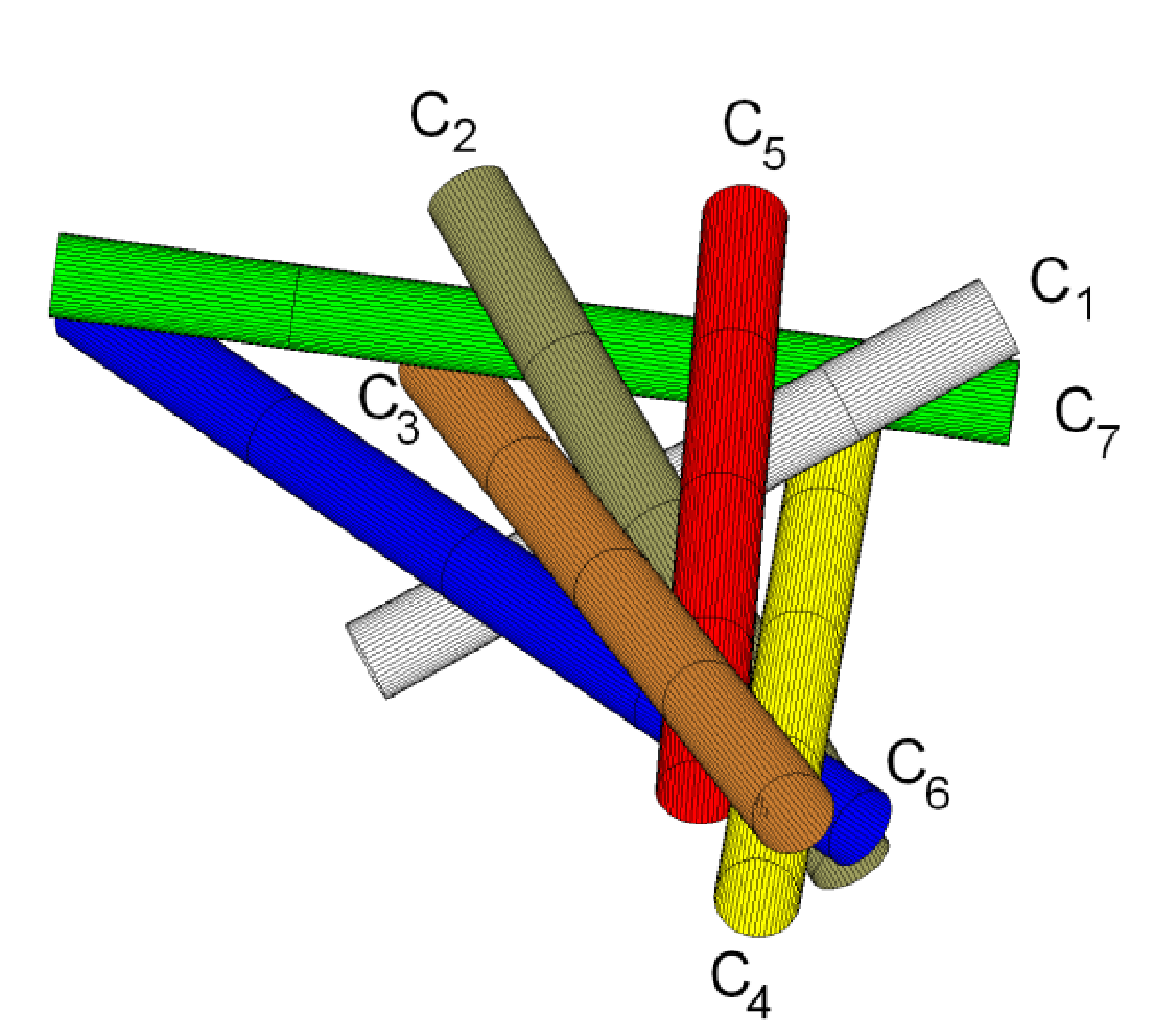}}}}
\put(-10,-40){\makebox{\textbf{Figure 2.} The second set of seven mutually touching inf{\kern0pt}inite cylinders }}
\end{picture}
\end{center}
\newpage

\section{Verif{\kern0pt}ication of the roots}

HOM4PS-2.0 provides the solution up to 50 digits
(the f{\kern0pt}irst 12 of which being correct),
that can be used as a starting point
of a solver using f{\kern0pt}loating-point arithmetic like
 \texttt{fsolve} in Maple 13.
With several accuracy levels adjusted previously by
\texttt{Digits:=$10^r$} ($r=2,3,4$),
CPU times of running \texttt{fsolve} without any further
specif{\kern0pt}ication on a personal computer with
Pentium(R) 4 CPU 3.4GHz and 2GB of RAM
are listed in Table 2.

\begin{center}
\begin{tabular}{|c|c|c|c|}
\hline
\verb"Digits" &      $10^2$          &          $10^3$        &        $10^4$       \\
\hline
    CPU time  &    0.6 seconds       &        4 seconds       &       130 seconds   \\
\hline
\end{tabular}
\end{center}
\begin{center}
\textbf{Table 2.} CPU time of \texttt{fsolve} for $10^r$ ($r=2,3,4$) correct digits
\end{center}

It's worth noting that \texttt{fsolve} recovers the solutions,
within approximately the same CPU time as in Table 2,
even if the starting values are truncated at 2 digits.
Moreover, truncation at 1 digit still works for the
f{\kern0pt}irst solution. Truncation at 1 digit, except for
$ t_3  $, which is truncated at 2 digits ($-0.03$),
works for the second solution.

However, a large number of correct digits is still
not mathematical correctness. Two exact verif{\kern0pt}ication
methods, alphaCertif{\kern0pt}ied and the interval Krawczyk method
are applied.
Any of them would be suf{\kern0pt}f{\kern0pt}icient of its own, nevertheless,
two is at least not worse than one.

\subsection{alphaCertif{\kern0pt}ied}

Smale's $\alpha$-theory  \cite{Smale1986} provides a positive, ef{\kern0pt}fectively computable
constant $\alpha(F, \mathbf{x})$ for a polynomial system
$F:\mathbb{C}^n \rightarrow \mathbb{C}^n$ and a point
$\mathbf{x}\in \mathbb{C}^n$ with the property that if
\[
\alpha(F, \mathbf{x})\leq \frac{13-3\sqrt{17}}{4}\approx 0.1576,
\]
then Newton's iteration starting from $\mathbf{x}$ converges quadratically
to a solution $\boldsymbol\xi$  close to $\mathbf{x}$ of the system $
F =\mathbf{0}$.
Based on Smale's theory
Hauenstein and Sottile \cite{HauensteinSottile2012} developed algorithms which, for given
$F$ and $\mathbf{x}$, compute an upper bound on $\alpha(F,\mathbf{x})$
and on some related quantities. On that basis they have
built a multipurpose verif{\kern0pt}ication software called alphaCertif{\kern0pt}ied. It can
produce a certif{\kern0pt}icate that \\
(i) $\mathbf{x}$ is an approximate solution
of $F=\mathbf{0}$ in the above sense;  \\
(ii) an approximate solution corresponds to an isolated solution; \\
(iii) the solution $\boldsymbol\xi$ corresponding to $\mathbf{x}$ is real (for real
$F$).

We have used alphaCertif{\kern0pt}ied v1.2.0 (August 15, 2011,
GMP v4.3.1 \& MPFR v2.4.1-p5) with Maple 13 interface.
The input of alphaCertif{\kern0pt}ied is system (\ref{eq:line1-linej})-(\ref{eq:linei-linej})
and the approximate solutions in Table 1.
We need to write the f{\kern0pt}irst solution up to at least 12 digits,
otherwise algorithm alphaCertif{\kern0pt}ied does not certify it.
The output of alphaCertif{\kern0pt}ied with the f{\kern0pt}irst solution
as in Table 1 consists of
$\alpha = 4.4333\cdot10^{-2},$
$\beta =  3.1668\cdot10^{-12},$
$\gamma = 1.3999\cdot10^{10}$
(see \cite{HauensteinSottile2012} for the details of
$\alpha, \beta, \gamma$).
The second solution has to be written up to at least
11 digits in order to be certif{\kern0pt}ied.
The output of alphaCertif{\kern0pt}ied with the second solution
(truncated at 11 digits) consists of
$\alpha = 6.578\cdot10^{-2},$
$\beta =  2.2387\cdot10^{-11},$
$\gamma = 2.9392\cdot10^{9}.$
Both solutions have been certif{\kern0pt}ied to be real and isolated solutions.

\subsection{The interval Krawczyk method}

We seek for real solutions among the numerical solutions with imaginary parts
less than the heuristic threshold $\theta =10^{-8}$.
The residuals of the real solutions are less than $5\cdot10^{-14}$ and their condition
numbers are at most $4.8\cdot10^{4}$, which show that these solutions are
numerically reliable.

To guarantee that in a small neighborhood of each numerical solution there
is a unique exact physical solution, the interval Krawczyk method \cite{Krawczyk1969}
is applied for verif{\kern0pt}ication. The method is based on the
following fact:~ for
a smooth function $F:\mathbb{R} ^{n}\rightarrow \mathbb{R}^{n}$
and a point $\mathbf{x}\in \mathbb{R}^{n}$, let ~$[\;\mathbf{x}%
\;]_{r}\subset \mathbb{R}^{n}$~ be the ball centered at $\mathbf{x}$
with radius $r>0$.~ Namely,
\begin{equation*}
[\;\mathbf{x}\;]_{r}~=~\left\{ \mathbf{y}\in \mathbb{R}^{n}:\left\| \mathbf{y%
}-\mathbf{x}\right\| _{\infty }\leq r\right\} ,
\end{equation*}
where $\left\| \mathbf{ \,\, }\right\| _{\infty }$ is the inf{\kern0pt}inity norm.
Assuming that the derivative of $F$ at $\mathbf{x}$, denoted by ~$DF(\mathbf{x})$,
~is nonsingular, the Krawczyk set of $F$ associated with ~$[\;\mathbf{x}%
\;]_{r}$~ is def{\kern0pt}ined as
\begin{equation*}
K(F,[\;\mathbf{x}\;]_{r})~=~\mathbf{x}-DF(\mathbf{x})^{-1}F(\mathbf{x})+%
\left[ I-DF(\mathbf{x})^{-1}DF([\;\mathbf{x}\;]_{r})\right] ([\;\mathbf{x}%
\;]_{r}-\mathbf{x}).
\end{equation*}
If the Krawczyk set is contained in the interior of ~$[\;\mathbf{x}\;]_{r}$%
,~ then there exists a unique zero of $F$ in ~$[\;\mathbf{x}\;]_{r}$.

The task of verif{\kern0pt}ication is implemented by using the interval arithmetic in
INTLAB (INTerval LABoratory) \cite{Rump-INTLAB}. ~In this implementation each
numerical solution $\mathbf{x}$ is taken as the center of the ball
~$[\;\mathbf{x}\;]_{r}$~ with radius $r=10^{-8}$.
Again, both solutions have been certif{\kern0pt}ied to be real and isolated solutions.

\section{Conclusions and open questions}

It remains an open question
whether seven is the maximal number of mutually touching
inf{\kern0pt}inite cylinders. Following the same idea for eight cylinders,
a polynomial system of 25 variables and 27 equations is resulted in.
It is not yet dis/proven whether it has a solution.
In case of seven cylinders, alternative choices instead of that
the f{\kern0pt}irst two cylinders are orthogonal need to be analyzed.
The maximal number of lines in $\mathbb{R}^n$
($n  > 3$) having the same pairwise nonzero distance is also unknown.
The authors believe that the method proposed can be applied for a
wide class of similar geometrical problems.

The angle of the f{\kern0pt}irst two cylinders was f{\kern0pt}ixed at $90^{\circ}$
in (\ref{eq:specificationline1})-(\ref{eq:specificationline2})
in order to have the same number of variables and equations.
A natural question arises whether the solutions we have found
can be extended if the angle of the f{\kern0pt}irst two cylinders varies.
We hope to return to this problem with an af{\kern0pt}f{\kern0pt}irmative answer.

Peter V.~Pikhitsa (Seoul National University) contacted
us after we had uploaded our manuscript to arXiv. He let us know
how physicists, including himself, investigated the problem
of mutually touching inf{\kern0pt}inite cylinders
\cite{Pikhitsa2004,Pikhitsa2007,PikhitsaChoiKimAhn2009,PikhitsaChoi2014}.
We think that the approximate solutions they found, especially the one
in \cite[Fig.1b]{PikhitsaChoiKimAhn2009} could be used
as starting points of the ref{\kern0pt}inement process,
which may lead to exact solutions, similar to ours.
The generalized problem that allows dif{\kern0pt}ferent diameters of the
cylinders is of mathematical interest, too.

\section{Acknowledgements}
The authors are grateful to Andr\'as Bezdek (Auburn University;
Alfr\'ed R\'enyi Institute of Mathematics, Hungarian Academy of Sciences)
for communicating on existing results related to Littlewood's
problem, and to \'Agnes Sz\'ant\'o (NC State University)
for useful discussion on $\alpha$-theory.
We thank Peter V.~Pikhitsa (Seoul National University) for
informing us about his related work.
We also thank the anonymous reviewer for valuable suggestions and comments.

\bibliographystyle{elsarticle-num}
\bibliography{<your-bib-database>}

\begin{thebibliography}{99}



\bibitem{AmbrusBezdek2008} 
G. Ambrus, A. Bezdek, On the number of mutually touching cylinders. Is it
8?, European Journal of Combinatorics 29 (8) (2008) 1803\textendash1807.

\bibitem{Bezdek2005} A.~Bezdek, On the number of mutually touching
cylinders, Combinatorial and Computational Geometry, MSRI Publication, 52
(2005) 121\textendash127.

\bibitem{BrassMoserPach2005}
P.~Brass, W.~Moser, J.~Pach,
Research Problems in Discrete Geometry,
Springer, 2005.

\bibitem{ChenLeeLi2014}
T.~Chen, T.L.~Lee, T.Y.~Li,
Mixed Volume Computation in Parallel,
Taiwanese Journal of Mathematics
18(1)
(2014)
93\textendash114.

\bibitem{Gardner1988} 
M.~Gardner, Hexaf{\kern0pt}lexagons and Other Mathematical Diversions: The
First Scientif{\kern0pt}ic American Book of Puzzles and Games, University of
Chicago Press, 1988, pp.~110\textendash115.

\bibitem{GreenspanBenney1997} 
H.P.~Greenspan, D.J.~Benney, Calculus: An Introduction to Applied
Mathematics, Breukelen Press, Brookline, MA, 1997, pp.~432\textendash433.

\bibitem{HauensteinSottile2012} 
J.D.~Hauenstein, F.~Sottile, Algorithm 921: alphaCertif{\kern0pt}ied:
Certifying Solutions to Polynomial Systems, ACM Transactions on Mathematical
Software, 38 (4) (2012) Article 28. DOI 10.1145/2331130.2331136

\bibitem{HuberSturmfels1995}
B.~Huber, B.~Sturmfels,
A polyhedral method for solving sparse polynomial systems,
 Mathematics of Computation
64
 (212)
 (1995)
 1541\textendash1555.

\bibitem{Krawczyk1969}
R.~Krawczyk,
Newton-Algorithmen zur Bestimmung von Nullstellen mit Fehlerschranken,
Computing 4 (3) (1969) 187\textendash201.

\bibitem{LeeLiTsai2008}
T.L.~Lee, T.Y.~Li, C.H.~Tsai, HOM4PS-2.0, A software package for solving
polynomial systems by the polyhedral homotopy continuation method, Computing
83 (2008) 109\textendash133. 

\bibitem{Littlewood1968} J.E.~Littlewood, Some problems in real and complex
analysis, Heath Mathematical Monographs, Raytheon Education, Lexington,
Massachusetts, 1968. 

\bibitem{Ogilvy1962}
C.S.~Ogilvy,
Tomorrow's math: unsolved problems for the amateur,
Oxford University Press, New York, 1962, pp.~60-61 and p.~153.

\bibitem{Pikhitsa2004}
P.V.~Pikhitsa,
Regular network of contacting cylinders with implications for materials
with negative Poisson ratios,
Physical Review Letters 93 (1) (2004) Article 015505.
DOI 10.1103/PhysRevLett.93.015505

\bibitem{Pikhitsa2007}
P.V.~Pikhitsa,
Architecture of cylinders with implications for materials with negative Poisson ratio,
Physica Status Solidi B 244 (3) (2007) 1004\textendash1007.

\bibitem{PikhitsaChoiKimAhn2009}
P.V.~Pikhitsa, M.~Choi, H-J.~Kim, S-H.~Ahn,
Auxetic lattice of multipods,
Physica Status Solidi B 246 (9) (2009) 2098\textendash2101.

\bibitem{PikhitsaChoi2014}
P.V.~Pikhitsa, M.~Choi,
Seven, eight, and nine mutually touching inf{\kern0pt}initely
long straight round cylinders: Entanglement in
Euclidean space,
manuscript, arXiv:1312.6207

\bibitem{Rump1999} 
S.M.~Rump, INTLAB -- INTerval LABoratory, in: T.~Csendes, editor,
Developments in Reliable Computing, Kluwer Academic Publishers, Dordrecht,
1999, pp.~77\textendash104.

\bibitem{Rump-INTLAB} S.M.~Rump, INTLAB -- INTerval LABoratory. \newline
{\small {\ \verb|http://www.ti3.tu-harburg.de/rump/intlab/|} \newline
}{\normalsize Accessed on 22 August 2013 }

\bibitem{Smale1986}
\normalsize
S.~Smale, Newton's method estimates from data at one point, in
R.E.~Ewing, K.I.~Gross, C.F.~Martin (editors): The Merging of Disciplines:
New Directions in Pure, Applied, and Computational Mathematics, Springer,
New York, 1986, pp.~185\textendash196. 

\bibitem{Weisstein-MathWorld} {\normalsize 
E.W.~Weisstein, ``Line-Line Distance.'' From MathWorld--A Wolfram Web
Resource. \newline
}{\small {\ \verb|http://mathworld.wolfram.com/Line-LineDistance.html|}
\newline
}{\normalsize Accessed on 22 August 2013 }

\end{thebibliography}







\end{document}